\documentclass{article}

\usepackage{amsmath}
\usepackage{palatino}
\usepackage{amssymb}
\usepackage{mathrsfs}
\usepackage[all]{xy}
\usepackage{graphicx}
\usepackage{latexsym}
\usepackage{amsthm}
\usepackage{amscd}
\usepackage{mathrsfs}

\newtheorem{proposition}{Proposition}
\newtheorem{theorem}[proposition]{Theorem}
\newtheorem{lemma}[proposition]{Lemma}
\newtheorem{corollary}[proposition]{Corollary}
\theoremstyle{definition}
\newtheorem{definition}[proposition]{Definition}
\newtheorem{remark}[proposition]{Remark}

\def\q{{\mathbf Q}}

\def\qbar{\bar{\mathbf Q}}

\DeclareMathOperator\res{Res}

\title{On the Belyi degree of a number field}
\author{Leonardo Zapponi}
\date{\today}

\begin{document}

\maketitle
\begin{abstract} In this short note we introduce the Belyi degree of a number field $K$, which is the smallest degree of a dessin d'enfant having $K$ as field of moduli. After the description of some general properties (for example, the fact that there exist finitely many number fields of bounded Belyi degree), we give a lower and an upper bound for such an invariant. We finally give some explicit examples for quadratic fields.
\end{abstract}

\section*{Introduction}

Let $K$ be a number field. By a well-known result there exists a dessin d'enfant (cover of the projective line unramified outside three points) having $K$ as field of moduli. The Belyi degree of $K$ is the minimal degree of such a cover. The present note is an invitation to the study of this invariant.

The paper is organized as follows:

In the first section, we give some general results concerning the Belyi degree, as for example the fact that there exist finitely many number fields of bounded (Belyi) degree. 

In \S2, we obtain a lower bound for such an invariant (Theorem 4), which follows from a result of S. Beckmannn concerning the ramified primes in the field of moduli of a cover of the projective line defined over a number field. Although this bound seems rough, it can be shown that it is 'optimal'. Indeed, in \S4 we prove that it is reached for infinitely many imaginary quadratic fields.

In \S3 we use a construction of R. Li\c tcanu and a result of D. Roy and J. L. Thunder in order to obtain an upper bound of the Belyi degree (Theorem 7), only depending on the degree and on the (absolute value) of the discriminant of the number field. More precisely, we prove that once we have fixed the degree, the growth of the Belyi degree is at most polynomial on the discriminant. This leads to the introduction of the minimal discriminant exponent.

In \S4, which is a computational approach to the Belyi degree of quadratic fields, we obtain an upper bound which is linear on the discriminant and in some cases, we are also able to completely determine this invariant. As a consequence we prove that the minimal discriminant exponent is equal to $1$ for these fields.

We also include two appendices: in the  first we generalize two results of H. W. Lenstra which are very useful when studying  the field of moduli of a genus $0$ dessin d'enfant;  in the second we give the explicit equations for the dessins d'enfants introduced in \S4.

\section{General results}

Even if it is a well-known fact (see for example~\cite{Schneps}, or more generally~\cite{Dessin} for an introduction to dessins d'enfants) we start this section by proving that any number field appears as the field of moduli of a dessin d'enfant.

\begin{lemma} Let $K$ be a number field. Then there exists a dessin d'enfant having $K$ as field of moduli.
\end{lemma}

\begin{proof} By the Primitive Element Theorem, there exists an element $\alpha\in\bar{\mathbf Q}$ such that $K=\mathbf Q(\alpha)$. Let $E$ be an elliptic curve having $\alpha$ as $j$-invariant.  Following~\cite{Silverman}, $K$ is the smallest field of definition of $E$. Fix a model of $E$ defined over $K$ and a rational function $f\in K(E)$, which induces a cover $E\to\mathbf P^1$ defined over $K$. Following the algorithm in~\cite{Belyi}, there exists a polynomial $g\in\mathbf Q[X]$ such that the composite cover $h=g\circ f$ induces a cover of the projective  line unramified outside $\infty,0$ and $1$. Consider now an element $\sigma\in\mbox{Gal}(\bar{\mathbf Q}/\mathbf Q)$. Then $\sigma$ belongs to $\mbox{Gal}(\bar{\mathbf Q}/\mathbf K)$ if and only if $\sigma(\alpha)=\alpha$. In this case, by construction (since $E$ and $f$ are defined over $K$ and $g$ is defined over $\mathbf Q$), we find $^\sigma h=h$ and thus the cover $h$ is defined over $K$. If $\sigma(\alpha)\neq\alpha$ then the conjugate cover $^\sigma h$ cannot  be  isomorphic to $h$, since $E$ and $^\sigma E$ have different $j$-invariants. This shows that $K$ is the field of moduli (and even the minimal field of definition) of $h$. 
\end{proof}

\begin{definition} The {\it Belyi degree} of a number field $K$ is the smallest degree of a dessin d'enfant having $K$ as field of moduli.
\end{definition}

\begin{lemma}\label{prop1} There exist finitely many number fields of fixed Belyi degree. Two isomorphic number fields have the same Belyi degree.
\end{lemma}

\begin{proof} The first assertion directly follows from the fact that there are finitely many dessins d'enfants of given degree. If $\sigma:K\to L$ is an isomorphism of number fields and $D$ is a dessin d'enfant having $K$ as field of moduli, then $^\sigma D$ has $L$ as field of moduli.
\end{proof}

\section{A lower bound}

\begin{theorem}\label{prop2} For any number field $K$, we have the inequality
$$\deg_B(K)\geq p,$$
where $p$ is the greatest prime which ramifies in $K$.
\end{theorem}

\begin{proof} Following a result of S. Beckmann~\cite{Beckmann}, the field of moduli $K$ of a dessin d'enfant is unramified outside the primes dividing the order of its monodromy group. In particular, if the degree of the dessin d'enfant is strictly less than a prime $p$, then the order of its monodromy group divides $d!$ and thus $p$ does not ramifies in $K$.
\end{proof}

This result seems quite rough but we will see in \S4 that this bound is reached for infinitely many (imaginary quadratic) number fields.

\section{An upper bound. The minimal discriminant exponent}

In this section, we follow the ideas and techniques of R. Li\c tcanu in~\cite{Litcanu} and then use a result of D. Roy and J. L. Thunder in~\cite{Roy}. We start with some preliminary results on heights.

Let $K$ be a number field and consider an element $x=[x_0,\dots,x_n]\in\mathbf P^n(K)$. Recall that the (multiplicative) height of $x$ is defined as
$$H_K(x)=\prod_{v}\max_i\{|x_i|_v\},$$
where the product is taken over all the places of $K$ (in order to avoid exponents, we distinguish two conjugated places). The normalized height is defined as
$$H(x)=H_K(x)^{1/d},$$
where $d=[K:\mathbf Q]$. For an element $x\in K$ we set $H(x)=H([x,1])$. Recall also that the height of a polynomial $f=a_0+\cdots+a_nX^n\in\mathbf K[X]$ is defined as $H(f)=H([a_0,\dots,a_n])$. In particular if $a_0,\dots,a_n$ are integers with no common divisor, we simply get
$$H(f)=\max_i\{|a_i|\}.$$

\begin{lemma}\label{lemma2} For any polynomial $f=a_0+\cdots+a_nX^n\in K[X]$, we have the inequality
$$H(f')\leq nH(f).$$
\end{lemma}

\begin{proof} We clearly have the identity
$$H(f')=H([a_1,2a_2,\dots,na_n]).$$
If $v$ is a finite place of $K$, we obtain
$$\max_i\{|ia_i|_v\}=\max_i\{|i|_v|a_i|_v\}\leq\max_i\{|a_i|_v\},$$
while for $v$ infinite we get
$$\max_i\{|ia_i|_v\}\leq\max_i\{|na_i|_v\}=n\max_i\{|a_i|_v\}$$
and the result follows by taking the product and then $d$-th roots.
\end{proof}

We will also use the fundamental inequalities
\begin{equation}\label{eq1}
H(f(x))\leq(n+1)H(f)H(x)^n,
\end{equation}
which holds for any $f\in K[X]$ of degree $n$ and any $x\in K$ (see for example~\cite{Serre}, p. 13) and
\begin{equation}\label{eq2}
2^{-n}\prod_{i=1}^nH(x_i)\leq H(f)\leq2^{n-1}\prod_{i=1}^nH(x_i),
\end{equation}
for any polynomial $f=c\prod_{i=1}^n(X-x_i)\in K[X]$ (cf.~\cite{Silverman}, Theorem 5.9).

\begin{lemma}\label{lemma3} Consider a monic polynomial $f\in K[X]$ of degree $n$ and the factorization $f'=n\prod_i(X-y_i)$. Setting
$$\hat f=\prod_i(X-f(y_i))=(-1)^{n-1}n^{-n}\res_Y(X-f(Y),f'(Y))\in K[X],$$
we have $\deg(\hat f)=n-1$ and the inequality
$$H(\hat f)\leq2^{n^2}(n+1)^{2n}H(f)^{2n}.$$
\end{lemma}

\begin{proof} Following Lemma~\ref{lemma2} and the above inequalities, we obtain
$$\aligned 
H(\hat f)&\leq2^{n-2}\prod_iH(f(y_i))\leq2^{n-2}(n+1)^{n-1}H(f)^{n-1}\prod_iH(y_i)^n\leq\\
&\leq2^{n^2-2}(n+1)^{n-1}H(f) ^{n-1}H(f')^n\leq2^{n^2-2}n^n(n+1)^{n-1}H(f)^{2n-1}
\endaligned$$
from which the lemma follows.
\end{proof}

\begin{theorem}\label{upper_bound} The Belyi degree of a number field $K$ is bounded by a constant only depending on its degree $n$ and on its discriminant $\Delta$. More precisely, there exist two effective constants $a,b\in\mathbf R$ only depending on $n$ such that we have the inequality
$$\deg_B(K)\leq a|\Delta|^b.$$
\end{theorem}

\begin{proof} Let $1,x_2,\dots,x_n$ be a basis of $K$ over $\mathbf Q$ and consider the polynomial
$$f_0=X+X^2+x_2X^3+\cdots+x_nX^{n+1}+X^{n+3}\in K[X].$$
The cover induced by $f_0$ is unramified outside the set $S_0=\{\infty,f_0(y_1),\dots,f_0(y_{n+2})\}$, where $y_1,\dots,y_{n+2}$ are the roots of $f_0'$. Following the notation of Lemma~\ref{lemma3}, consider the polynomial
$$f_1=N_{K/\q}(\hat f_0)\in\q[X]$$
of degree $n(n+2)<(n+1)^ 2$ and for any $m\in\{2,\dots,n(n+2)\}$ define $f_m\in\q[X]$ recursively by the relation $f_{m+1}=\hat f_m$. Set also
$$S_m=f_m(S_{m-1}\cup\{z_i\,\,|\,\, f_m'(z_i)=0\})=f_m(S_{m-1})\cup\{z_i\,\,|\,\, f_{m+1}(z_i)=0\},$$
so that $S_m$ has cardinality less than or equal to $(n+1)^2+1$ and the cover induced by the polynomial $g_m=f_m\circ\cdots f_1\circ f_0$ is unramified outside this set. Since
$$\deg(f_m)=n(n+2)-m+1<(n+1)^2-1,$$
following Lemma~\ref{lemma3}, we find
$$\aligned
H(f_m)&\leq2^{(n+1)^4}(n+1)^{2(n+1)^2}H(f_{m-1})^{2(n+1)^2}\leq\\
&\leq\left[2^{(n+1)^4}(n+1)^{2(n+1)^2}\right]^{1+2(n+1)^2+\cdots+(2(n+1)^2)^{m-2}}H(f_1)^{2^{m-1}(n+1)^{2(m-1)}}\\
&\leq\left[2^{(n+1)^4}(n+1)^{2(n+1)^2}H(f_1)\right]^{2^{m-1}(n+1)^{2(m-1)}}=a_0H(f_1)^{b_0}
\endaligned$$
with $a_0$ and $b_0$ only depending on $n$. Since
$$H(f_1)\leq4^{(n+1)^2}H(\hat f_0)^n\leq4^{(n+1)^2}\left(2^{(n+3)^2}(n+4)^{2n+6}H(f_0)^{2n+6}\right)^n=a_1H(f_0)^{b_1},$$
where $a_1$ and $b_1$ only depend on $n$, we then obtain the inequality
$$H(f_m)\leq a_0(a_1H(f_0)^{b_1})^{b_0}=a_2H(f_0)^{b_2},$$
where, once again, the constants $a_2$ and $b_2$ only depend on $n$. By construction, the set $S_{n(n+2)}$ is contained in $\mathbf P^1(\q)$ and we have
$$\aligned
H(S_{n(n+2)})&=\max_{x\in S_{n(n+2)}}H(x)\leq H(f_{n(n+2)}\circ\cdots\circ f_2(0))\leq\\
&\leq H(f_{n(n+2)})H(f_{n(n+2)-1})\cdots H(f_2)^{(n^2+2n-2)!}\leq a_3H(f_0) ^{b_3}
\endaligned$$
with $a_3,b_3$ only depending on $n$. Now, following~\cite{Litcanu}, Lemme 4.1 (and its proof), there exists a rational function $h\in\q(X)$ such that the cover induced by the rational function $f=h\circ g_{n(n+2)}$ is unramified outside $\infty,0$ and $1$. Moreover, we have the bound
$$\deg(f)=\deg(h)\deg(g_{n(n+2)})\leq(n+1)^2!\max_{x\in S_{n(n+2)}}\{H(x)\}^{(n+1)^2(n^2+2n+4)}\leq a_4H(f_0)^{b_4}.$$
Now, following~\cite{Roy}, Theorem 2, there exist $f_0$ such that
$$H(f_0)\leq2^{3d+1}|\Delta|^{1/d},$$
so that we obtain the desired inequality.
We just have to check that the field of moduli of the cover induced by $f$ is $K$. This directly follows from Appendix A: by construction, $f$ is defined over $K$. Let $\sigma\in\mbox{Gal}(\qbar/\q)$ not belonging to $\mbox{Gal}(\qbar/K)$. We must prove that the cover induced by $^\sigma f$ is not isomorphic to the one induced by $f$. Suppose that the opposite holds. This can be restated as
$$^\sigma f=f\circ\varphi$$
with $\varphi\in\mbox{Aut}(\mathbf P^1)$. Remark that, by construction, the point at infinity is the only ramified point of $f$ in its fiber (above $1$). This means that $\varphi(\infty)=\infty$ and thus $\varphi=uX+b$ with $u\in\qbar^\times$ and $v\in\qbar$. Setting $g=h\circ f_{n(n+2)}\circ\cdots\circ f_1$, we have
$$g\circ{^\sigma f_0}={^\sigma f}=f(uX+v)=g\circ f_0(uX+v)$$
and Lemma~\ref{Lenstra_2} implies that there exist $w,z\in\qbar$ such that
$$^\sigma f_0=wf_0(uX+v)+z.$$
Now, since the coefficient of $X^{n+2}$ vanishes for both $f_0$ and $^\sigma f_0$, we obtain $v=0$. Similarly, since $f_0(0)={^\sigma}f_0(0)=0$, we find $z=0$. Finally, comparing the coefficients of $X$ and $X^2$, we obtain
$$\left\{\aligned
&wu=1\\
&wu^2=1
\endaligned\right.$$
which gives $u=w=1$ and thus $^\sigma f_0=f_0$, which is impossible. This concludes the proof of the theorem.
\end{proof}

\begin{remark}  As it is done in~\cite{Litcanu}, it is possible to give an explicit bound for the constants $a$ and $b$. In order to obtain a 'light' proof, we decided to omit this computation.
\end{remark}

The above  theorem asserts that, once we have fixed the degree of the number field, the growth of its Belyi degree is at most polynomial on its discriminant. It is then natural to ask which is the least exponent for which such a inequality holds. This leads to the notion of the minimal discriminant exponent.

\begin{definition} Let $n>1$ be an integer. The {\it minimal discriminant exponent} is the real number $\delta(n)=\inf V_n$, where we have set
$$V_n=\left\{b\,\,|\,\,\exists\,\,a\mbox{ such that }\deg_B(K)\leq a|\Delta_K|^b\mbox{ for any }K\mbox{ with }[K:\q]=n\right\}.$$
We say that $\delta(n)$ is {\it effective} if it belongs to $V_n$.
\end{definition}

Remark that $\delta(n)$ is effective if and only if $V_n$ is a closed set (this is the usual problem of passing from $\delta(n)+\varepsilon$ do $\delta(n)$). This is also equivalent to the existence of a global bound for the constants $a$ appearing in the definition of $V_n$. As we will se in the next section, it is possible to determine the minimal discriminant exponent for quadratic fields.

\section{Examples}\label{examples}

\subsection{Quadratic imaginary fields} We start this list of example with the study of quadratic imaginary fields. As we will see, in some cases, it is possible to completely determine the Belyi degree.

\begin{proposition}\label{bound1} Let $d$ be a positive integer and set $K=\mathbf Q(\sqrt{-d})$. We then have the inequality
$$\deg_B(K)\leq 2d+4.$$
\end{proposition}

\begin{proof} The two dessins d'enfants in Figure 1 are conjugated and have $K$ as field of moduli; their degree being $2d+4$, the result follows.
\end{proof}

\vskip.4cm
\begin{center}
\includegraphics[scale=.4]{BelyiDegreeFigure1.epsf}\\
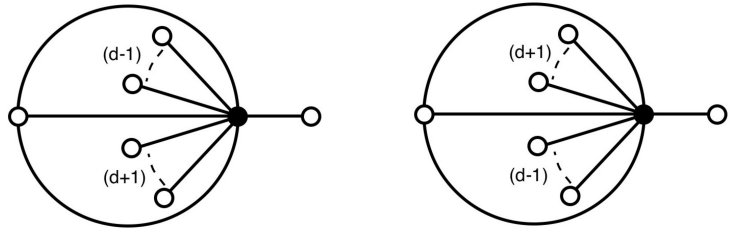
{Figure 1}
\end{center}
\vskip.4cm

The following result shows that the bound in Theorem~\ref{prop2} is sharp.

\begin{proposition} Let $a<b<c$ three positive integers such that $p=a+b+c$ is prime. Setting $K=\mathbf Q(\sqrt{-abcp})$, we have the identity
$$\deg_B(K)=p.$$
In particular, there exist infinitely many (quadratic imaginary) number fields such that the inequality in Theorem~\ref{prop2} is an equality.
\end{proposition}

\begin{proof} Since $p$ does not divide $abc$, it is the greatest prime which ramifies in $K$. In particular, Theorem~\ref{prop2} asserts that the Belyi degree of $K$ is greater than or equal to $p$. We just have to construct a dessin d'enfants of degree $p$ having $K$ as field of moduli. The two dessins d'enfants described in figure 2 are conjugated (as soon as $a,b$ and $c$ are pairwise distinct) of degree $p$ and with $K$ as field of moduli. The last part of the proposition simply follows from the fact that there exist infinitely many primes.
\end{proof}

\begin{corollary} Let $p>7$ be a prime not congruent to $1$ modulo $12$ and set $K=\mathbf Q(\sqrt{-p})$. We then have the identity
$$\deg_B(K)=p.$$
\end{corollary}

\begin{proof} The congruence condition implies that $p$ is congruent either to $2$ modulo $3$ or to $3$ modulo $4$. In the first case, let $p=3n+2$ and consider the dessin d'enfant described in the proof of the above proposition with $a=2, b=n$ and $c=2n$. In the second case, let $p=4n+3$ and take $a=3, b=n$ and $c=3n$.
\end{proof}

\vskip.4cm
\begin{center}
\includegraphics[scale=.4]{BelyiDegreeFigure2.epsf}\\
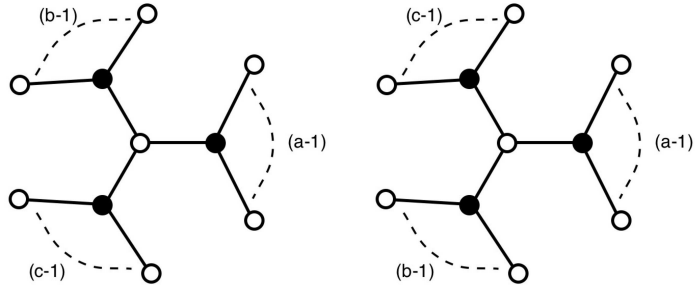
{Figure 2}
\end{center}
\vskip.4cm

\subsection{Real quadratic fields} We now study the case of real quadratic fields. The following result is similar to Proposition~\ref{bound1}, and even a little bit sharper.

\begin{proposition} Let $d\geq5$ be an integer and set $K=\mathbf Q(\sqrt{d})$. We then have the inequality
$$\deg_B(K)\leq 2d-2.$$
\end{proposition}

\begin{proof}The two dessins d'enfants in Figure 3 are conjugated and have $K$ as field of moduli; their degree being $2d-2$, the result follows.
\end{proof}

\vskip.4cm
\begin{center}
\includegraphics[scale=.35]{BelyiDegreeFigure3.epsf}\\
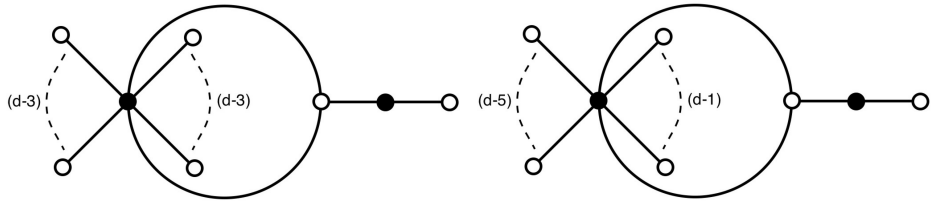
{Figure 3}
\end{center}
\vskip.4cm

\subsection{Computing the minimal discriminant exponent for quardatic fields}

The results of the last two paragraphs allow us to completely determine the minimal discriminant exponent for quadratic fields.

\begin{proposition} The minimal discriminant exponent for quadratic fields is equal to $\delta(2)=1$ and is effective.
\end{proposition}

\begin{proof} Set $K=\q(\sqrt{d})$ with $d$ squarefree, so that $|\Delta_K|\leq4|d|$. The last two bounds for the Belyi degree lead to
$$\deg_B(K)\leq24|\Delta_K|.$$
This implies that $\delta(2)\leq1\in V_2$ (cf. the notation of \S3). Suppose that $\delta(2)<1$. This means that there exists $\varepsilon>0$ and a constant $a$ such that
$$\deg_B(K)\leq a|\Delta_K|^{1-\varepsilon}.$$
Now, there exist infinitely many prime numbers $p$ not congruent to $1$ modulo $12$ and, setting $d=-p$, Corollary 12 implies that $\deg_B(K)=p\geq\frac14|\Delta_K|$. In particular, for large enough $p$, we find
$$|\Delta_K|^\varepsilon>4a,$$
and thus
$$\deg_B(K)\geq\frac14|\Delta_K|=\frac14|\Delta_K|^{1-\varepsilon}\cdot|\Delta_K|^\varepsilon>a|\Delta_K|^{1-\varepsilon},$$
which is impossible. This concludes the proof of the proposition.
\end{proof}

\section*{Appendix A. Two useful lemmas}

In this first appendix, we generalize two results of H. W. Lenstra in~\cite{Schneps}.

\begin{lemma}\label{Lenstra_1} Let $f\in K(X)$ be a rational function of degree $n=md$ and suppose that there exists a monic polynomial $h\in XK[X]$ of degree $d$ such that $f=g\circ h$ for a rational function $g\in K(X)$. Then $h$ is unique.
\end{lemma}

\begin{proof} Set $f=f_1/f_2$ with $f_1,f_2\in K[X]$ coprime and $g=g_1/g_2$ with $g_1,g_2\in K[X]$ coprime. We can suppose that $\deg(f_1)=n$. Since $g_1\circ h$ and $g_2\circ h$ are coprime, the identity
$$f_1\cdot g_2\circ h=f_2\cdot g_1\circ h$$
implies that there exists $u\in K^\times$ such that
$$f_1=ug_1\circ h=g_3\circ h,$$
where we have set $g_3=ug_1\in K[X]$. But in this case, following~\cite{Schneps}, Lemma II.2, $h$ is unique. 
\end{proof}

\begin{lemma}\label{Lenstra_2} Let $f_1,f_2\in K(X)$ and $h_1,h_2\in K[X]$ of the same degree such that $f_1\circ h_1=f_2\circ h_2$. Then there exist $a,b\in K$ such that $h_2=ah_1+b$.
\end{lemma}

\begin{proof} Let $a_i$ be the leading coefficient of $h_i$ and consider the monic polynomial
$$\tilde h_i=a_i^{-1}(h_i-h_i(0))\in XK[X].$$
In this case, setting
$$g_i=f_i(a_iX+h_i(0))\in K(X),$$ we obtain the identity
$$g_1\circ\tilde h_1=g_2\circ\tilde h_2$$
and the above lemma implies that $\tilde h_1=\tilde h_2$, from which the result follows.
\end{proof}

\section*{Appendix B. Explicit equations}

We close this paper by giving the explicit equations for the dessins d'enfants used in the proofs in~\S\ref{examples}. Although they look quite different, they can be studied in the same way. Consider three integers $a,b$ and $c$ with $abc(b+c)\neq0$ and $a+b+c>0$. We also assume that $a,b,c,a+b+c$ are pairwise distinct. We want to describe all the rational functions
$$f=(1-X)^a(1-xX)^b(1-yX)^c\in\bar{\mathbf Q}(X)$$
with $xy(x-y)(x-1)(y-1)\neq0$ such that the logarithmic derivative $df/f$  has a unique zero at the origin. In this case, $f$ induces a cover $\mathbf P^1\to\mathbf P^1$ unramified outside $\infty,0$ and $1$. Remark that the  condition $a+b+c>0$ implies that $f$ has a pole at infinity. It easily turns out that the condition on the  logarithmic derivative of $f$ is  equivalent to the  two equations
$$\left\{\aligned
&a+bx+cy=0\\
&a+bx^2+cy^2=0
\endaligned\right.$$
which lead to
$$\left\{\aligned
&x=-\frac{ab-\sqrt{-abc(a+b+c)}}{b(b+c)}\\
&y=-\frac{ac+\sqrt{-abc(a+b+c)}}{c(b+c)}
\endaligned\right.$$
If $\Delta=-abc(a+b+c)$ is not a perfect square then we have  $xy(x-y)(x-1)(y-1)\neq0$, as desired. Clearly, the rational function $f$ is defined over the field $K=\mathbf Q(\sqrt{\Delta})$. It is then easily shown that $K$ is in fact its field of moduli (this comes from the fact that $a,b,c,a+b+c$ are pairwise distinct).

 \end{document}